\documentclass[12pt]{article}
\usepackage{graphicx}
\usepackage{color}
\usepackage{amsmath}
\usepackage{dsfont}
\usepackage{placeins}
\usepackage{amssymb}
\usepackage{wasysym}
\usepackage{abstract}
\usepackage{hyperref}
\usepackage{parskip}

\usepackage[margin=1in]{geometry}

\newcounter{probnum}
\allowdisplaybreaks

\definecolor{gray}{rgb}{0.5,0.5,0.5}
\definecolor{black}{rgb}{0,0,0}
\definecolor{white}{rgb}{1,1,1}
\definecolor{blue}{rgb}{0.5,0.5,1}

\definecolor{green}{rgb}{0.133,0.545,0.133}

\definecolor{yellow}{rgb}{1,0.549,0}

\definecolor{red}{rgb}{1,0.133,0.133}

\definecolor{purple}{rgb}{0.58,0,0.827}

\definecolor{backgcode}{rgb}{0.97,0.97,0.8}
\definecolor{Brown}{cmyk}{0,0.81,1,0.60}
\definecolor{OliveGreen}{cmyk}{0.64,0,0.95,0.40}
\definecolor{CadetBlue}{cmyk}{0.62,0.57,0.23,0}

\newcommand{\beqn}{\vspace{-0.25cm}\begin{eqnarray*}}
\newcommand{\eeqn}{\end{eqnarray*}}
\newcommand{\bneqn}{\vspace{-0.25cm}\begin{eqnarray}}
\newcommand{\eneqn}{\end{eqnarray}}
\newcommand{\benum}{\begin{enumerate}}
\newcommand{\eenum}{\end{enumerate}}

\newcommand{\parens}[1]{\left(#1\right)}

\newcommand{\bracks}[1]{\left[#1\right]}
\newcommand{\braces}[1]{\left\{#1\right\}}

\newcommand{\expe}[1]{\mathbb{E}\bracks{#1}}

\renewcommand{\exp}[1]{\mathrm{exp}\parens{#1}}

\renewcommand{\min}[1]{\text{min}\braces{#1}}

\title{Pick up sticks}
\author{{\it Larry} Shepp, {\it Doron} Zeilberger, {\it Cun-Hui} Zhang}

\begin{document}
\maketitle

\section{Introduction}

Break a ruler of length $n$ inches with indentations at each inchmark into $n$ sticks of length one by repeatedly picking up all the fragments 
which are still of length greater than one and throwing them down all at once until all the sticks are of length one. 
Assume that each stick breaks into exactly two sticks uniformly at the inchmarks. What is the expected number of throws needed to complete the $n-1$ breaks?

We show that the answer to this question is $(4.31107\dots) \cdot \ln n \, \cdot (1+ o(1))$
by, as mathematicians are fond of doing, {\it reducing it to an already solved problem}.
The problem that we reduce it to-or rather show that it is almost trivially equivalent to- is the problem of finding
the asymptotic average height of a {\it binary search tree}.
The latter problem was first solved by  Luc Devroye[D], with a later refinement by Bruce Reed[Re].

We also present a possibly new (and we hope elegant) proof that $(4.31107\dots) \cdot \ln n \, \cdot (1+ o(1))$ is
an {\it upper bound} , first proved in 1979 by Robson[Ro].
Our approach is to first prove
a simple generating function for the number of throws needed to isolate the $j$\underline{th} stick and 
then use Chernoff's inequality.
The constant $\alpha = 4.31107\dots$ is the root of $\alpha \ln(2e/\alpha) = 1$. 

The {\it worst} case, if you are extremely unlucky, requires $n-1$ throws. In this case, first the stick breaks into a one-inch piece
and an $(n-1)$-inch piece, then the next throw breaks the $(n-1)$-inch stick into a one-inch piece and an $(n-2)$-inch piece and so on.

The {\it best} case (if $n=2^k$) is $\log_2 n$. First it breaks into two equal pieces each of length
$n/2$ inches, then into four pieces each of length $n/4$ inches etc., and in $k=\log_2 n$ throws into $n$ one-inch pieces.

But if you do it {\it many} times, and $n$ is large, what is the {\it expected} number of throws?
As we have already stated, the answer is $(4.31107 \dots) \cdot (\ln n) \cdot (1+o(1))$,
so the ratio of the asymptotic average to the best case is only
$(4.31107 \dots) \cdot \ln 2 \cdot (1+o(1))= 2.98821+o(1)$.

Recall that a {\it full binary tree} is a ```family tree'' (but with single-parenthood), with a {\it root}
(Eve), and where every vertex either has two children (called the left-child and the right-child) or
no children at all. Childless vertices are called {\it leaves}. 
The {\it height} of a full binary tree is the length of a maximal path from the root to a leaf.

Each scenario of completely breaking an $n$-inch stick into $n$ one-inch pieces can be naturally associated with a full binary tree
with $n$ leaves. The first drop of the interval $[0,n]$ results either in $[0,1]$ and $[1,n]$, or
$[0,2]$ and $[2,n]$, \dots, or $[0,n-1]$ and $[n-1,n]$. Now each of the two pieces goes its own way,
generating its own full binary subtree. Of course, the ``number of throws needed to completely break the stick''
is the height of the corresponding full binary tree.

To see all the  $14$ breaking scenarios for a five-inch rule go to: \hfill\break
{\tt http://www.math.rutgers.edu/\~{}zeilberg/shepp/S5.html} \quad .

To see all the  $42$ breaking scenarios for a six-inch rule go to: \hfill\break
{\tt http://www.math.rutgers.edu/\~{}zeilberg/shepp/S6.html} \quad .

Now to each internal vertex (i.e. a vertex that is {\it not} a leaf) assign the label ``number of leaves in the subtree whose root it is''. It is obvious
that the {\it probability} of the breaking  scenario corresponding to any given full binary tree is the
product of $1/(i-1)$ over all the labels of non-leaves, since a stick of length $i$ may be broken
into two smaller pieces in $i-1$ equally likely ways.

So the {\it expectation} of the random variable ``number of throws needed to completely break an $n$-inch ruler''
is nothing but the {\it expected height} of all full binary trees with $n$ leaves, under the above
probability distribution. Let's call it $a(n)$.

Recall (or go to wikipedia) that any list of $n$ different numbers gives rise to
a {\it binary search tree} of $n$ vertices where the first record, let's call it $i_1$ is put at the root,
and its left (respectively right) subtree consists, recursively of the binary search tree of
the sublist (in the same order) consisting of all entries smaller  (respectively larger) than $i_1$.
If a list is empty, then it corresponds to the empty tree with no vertices.

For example, the binary search tree corresponding to the permutation $[4,3,9,1,5, 6, 2,11, 7,8, 10,12]$ can be seen here:  \hfill\break
{\tt http://www.math.rutgers.edu/\~{}zeilberg/shepp/BS12.html} \quad .

It is obvious that a binary search tree corresponds to a (not-necesssarily full) binary tree,
where every vertex may either have no children, or only a left-child, or only a right-child,
or both left- and right-children with $n$ vertices.

To see all the  $14$ binary search trees for the four records $\{1,2,3,4\}$  go to: \hfill\break
{\tt http://www.math.rutgers.edu/\~{}zeilberg/shepp/BS4.html} \quad .

To see all the  $42$ binary search trees for the four records $\{1,2,3,4, 5\}$  go to: \hfill\break
{\tt http://www.math.rutgers.edu/\~{}zeilberg/shepp/BS5.html} \quad .

If you pick a permutation of $\{1, \dots, n\}$ uniformly at random, it is clear that
the probability of it having a certain binary search tree is simply the 
product of $1/i(v)$, over all vertices $v$ of the tree,
where $i(v)$ is defined to be ``the number of vertices in the subtree rooted at $v$''.
For example, the probability of the tree shown in: \hfill\break
{\tt http://www.math.rutgers.edu/\~{}zeilberg/shepp/BS12.html} \quad .

is
$$
(1/12)(1/3)(1/2)(1/1)(1/8)(1/4)(1/3)(1/2)(1/1)(1/3)(1/1)(1/1) \quad .
$$

This is the probability assigned to a random binary search tree, and the average is defined according to
that probability distribution.

In a famous paper in theoretical computer science, 
Bruce Reed[Re] improved on a seminal paper 
of Luc Devroye[D] (who only had the highest-order asymptotics) by proving that 
if $b(n)$ is the average height of  binary tree with $n$ vertices under the above probability distribution, then
$$
b(n)= \alpha \ln n \, - \, \beta  \ln \ln  n + O(1) \quad,
$$
where $\alpha=4.31107\dots$ is the unique root in $[2, \infty)$ of the equation
$$
\alpha \ln ((2e)/\alpha)=1 \quad ,
$$
and $\beta=1.953 \dots$.

We are almost done! It is trivial to see (e.g. by induction) that a full binary tree with $n$ leaves has
$n-1$ {\it internal vertices} (i.e. non-leaves). There is a well-known, obvious, bijection between
binary trees with $n-1$ vertices to full binary trees with $n$ leaves. If a vertex has both children
leave it alone. If it only has a left-child, create for it a right-child that is a leaf.
If it only has a right-child, create for it a left-child that is a leaf. If it has no children (i.e. is a leaf
in the original tree), create for it both a left- and right-child, both leaves.
To go back is even easier! {\it Remove all the leaves}!

By removing all the leaves, the height gets reduced by one, 
(and all the labels as well, but the probability stays the same, due to the different definitions in both cases)
so we have the following relation between
the stick problem and the well-known average height of a random search tree problem:
$$
a(n)=b(n-1)+1 \quad .
$$
Of course, the asymptotics of $a(n)$ is the same as that of $b(n)$.

\section{ A Quick Proof of the Upper Bound}

We now give a quick argument that reproves the inequality
$a(n) \leq  \alpha \ln n \cdot (1+o(1))$, first proved by Robson[Ro],
thereby proving {\it half} of the Devroye result.

We hope, in the future, to extend our approach to proving the other half, namely $a(n) \geq  \alpha \ln n \cdot (1+o(1))$, 
thereby yielding a hopefully simpler proof, or at least an alternative one, of Devroye's result.

\subsection{The generating function for the time to isolate the $j$\underline{th} stick}

Let $\tau_{j,n}, j=1,\ldots,n$ denote the number of throws needed until the $j$\underline{th} inch stick is isolated and let
$f_{j,n}(\rho)$ denote the generating function of $\tau_{j,n}$,
\beqn
f_{j,n}(\rho) = \expe{\rho^{\tau_{j,n}}}.
\eeqn
It is clear that $f_{1,1}(\rho) = 1$, since $\tau_{1,1} = 0$, and that the recurrence,
\beqn
f_{1,n} = \frac{1}{n-1} \sum_{k=1}^{n-1} \rho f_{1,k}(\rho)
\eeqn
holds for the generating function of the number of throws needed to isolate the {\it first} stick. 
An easy induction shows that
\beqn
f_{1,n}(\rho) = \prod_{l=1}^{n-1} (1+\frac{\rho-1}{l}).
\eeqn
Since the number of throws to isolate the $j$\underline{th} stick is the sum of the number of throws needed until the
first break at $j-1$ {\it plus} the number until the first break at the $j$\underline{th} inchmark, and these numbers
are independent, so we see that $f_{j,n}(\rho) = f_{1,j-1}(\rho) f_{1,n-j}(\rho)$ which gives the general result that
\beqn
f_{j,n}(\rho) = \prod_{l=1}^{j-1}(1+\frac{\rho-1}{l}) \prod_{l^\prime = 1}^{n-j}(1+\frac{\rho-1}{l^\prime}).
\eeqn

\subsection{Expectation}

We now prove that the expected total number of throws needed to completely break the stick 
is no greater than $$\alpha \ln(n+1) - 1.$$ 
Since $\max_{j\le n}\tau_{j,n}$ is the number of throws needed, this assertion is equivalent to 
$$
\sum_{x=1}^{n-1} P\Big\{\max_{j\le n}\tau_{j,n} \ge x \Big\} \le \alpha\ln(n+1) - 1. 
$$

The explicit formula of the generating function allows us to apply the Chernoff bound 
$$
P\big(\tau_{j,n} > x \big) \le \rho^{-x} f_{j,n}(\rho)
$$ 
to each $\tau_{j,n}$. Let $m=(n+1)/2$.  
For $\rho\ge 2$, $(1+(\rho-1)/k)\le (1+1/k)^{\rho-1}$, so that 
$$
f(\rho)\le \left\{\prod_{l=1}^{j-1}\Big(1+\frac{1}{l}\Big) \prod_{l^\prime = 1}^{n-j}\Big(1+\frac{1}{l^\prime}\Big)\right\}^{\rho-1}
= \{j(n-j+1)\}^{\rho-1}\le m^{2\rho-2}. 
$$ 
Let $x > \alpha\ln m$, $y = x/\ln m$ and $\rho = y/2$. Since $\rho > \alpha/2 > 2$, 

\beqn
P\big(\max_{j\le n}\tau_{j,n} \ge x \big) 
&\le&  \sum_{j=1}^n \rho^{-x} f_{j,n}(\rho)
\cr &\le& 2\,\exp{\ln m - y \ln m \ln\rho + (2\rho-2)\ln m} 
\cr &=& 2\,\exp{(-1 + y\ln(2e/y))\ln m}. 
\eeqn
Since $(d/dt)((-1 + t\ln(2e/t)) = \ln(2/t) \le \ln(2/\alpha)$ for $\alpha\le t\le y$, 
$$
(-1 + y\ln(2e/y))\ln m \le (y-\alpha)\ln(2/\alpha) \ln m 
= (x - \alpha\ln m)\ln(2/\alpha). 
$$
It follows that $P\big(\max_{j\le n}\tau_{j,n} \ge x \big) \le 2(2/\alpha)^{x-\alpha\ln m}$. 
Let $u$ be the decimal part of $\alpha\ln m - \ln 2/\ln(2/\alpha)$. We find 
\beqn
\sum_{x=1}^{n-1}P\big(\max_{j\le n}\tau_{j,n} \ge x \big) 
& \le &  \sum_{x=1}^{n-1}\min{1,2(2/\alpha)^{x-\alpha\ln m}}
\cr &\le&  \alpha\ln m - \ln 2/\ln(2/\alpha) - u + \sum_{k=1}^\infty (2/\alpha)^{k-u}. 
\eeqn
Since the function $ - u + \sum_{k=1}^\infty (2/\alpha)^{k-u}$ is convex in $u\in [0,1]$ and takes a greater value at $u=1$ 
than $u=0$, its maximum, $1/(1-2/\alpha)-1=2/(\alpha-2)$, is attained at $u=1$. Consequently, the expectation 
of $\max_{j\le n}\tau_{j,n}$ is no greater than 
$$
\alpha\ln(n+1) - \alpha\ln 2 - \ln 2/\ln(2/\alpha) + 2/(\alpha-2). 
$$
The conclusion follows from $- \alpha\ln 2 - \ln 2/\ln(2/\alpha) + 2/(\alpha-2)< -1$. 

A slightly different simple argument using the Poisson dominance of Bernoulli yields 
an upper bound with the $\ln\ln n$ term: 
$$
E \max_{j\le n}\tau_{j,n} \le \alpha\ln  n - (\beta/3)\ln\ln n + O(1)
$$
with $\beta = (3/2)/\log(\alpha/2) = 1.953026$. However, the coefficient for the second order term,  
$\beta/3$, is not sharp according Reed's refinement of Devroye's result that states that

{\bf Remark}. The individual times, $\tau_j$, are rather smaller than the maximum of the $\tau$'s. Indeed, using the generating function of $\tau_j$
it is not hard to show that as $n \rightarrow \infty$,

\beqn
\tau_j = 2 \log{n} + \sqrt{2 \log{n}} \eta_j + o(\sqrt{log{n}})
\eeqn

where $\eta_j$ is standard normal random variable. Thus, it is only the relative independence of the $\tau$'s that allows the maximum to be as large as in Reed's theorem, $\sim 4.13\ldots \log{n}$.

{\bf Acknowledgement}: The stick problem was inspired by
the interesting paper [IK].

\section{Maple Packages}

This article is accompanied by two Maple packges for experimenting, and simulating, both
binary trees under the uniform distribution (not discussed in this article), and for the
stick problem (alias binary search trees). The most important one is
{\tt BinaryTrees} and the lesser one is {\tt ArbresBinaires}. 
There is yet another Maple package, {\tt EtzBinary}, mainly for plotting, and that's
how we drew the pictures linked to above. All three packages, as well as
some sample input and output files, can be gotten from the ``front'' of this article: \hfill\break

{\tt http://www.math.rutgers.edu/\~{}zeilberg/mamarim/mamarimhtml/stick.html} .

\section{ References}

[D] Luc Devroye, {\it A note on the height of a binary search tree}, J. of the ACM {\bf 33}(3) (July 1989), 489-493.

[IK] Yoshiaki Itoh and  P. L. Krapivsky,
{\it Continuum cascade model of directed random graphs: traveling wave analysis}, {\tt http://arxiv.org/abs/1206.3711 } .

[Re] Bruce Reed, {\it The height of a random binary search tree},
Journal of the ACM {\bf 50} (2003), 306-332.

[Ro] J.M. Robson, {\it The height of binary search trees}, Aust. Computing J. {\bf 11} (1979), 151-153.

\bigskip
\hrule
\bigskip
Lawrence Shepp, Wharton School, University of Pennsylvania,  {\tt shepp@wharton.upenn.edu}
\smallskip
Doron Zeilberger, Department of Mathematics, Rutgers University, {\tt zeilberg@math.rutgers.edu}
\smallskip
Cun-Hui Zhang, Department of Statistics, Rutgers University, {\tt czhang@stat.rutgers.edu}
\end{document}